\documentclass[12pt]{article}
\usepackage{amsmath,amssymb,theorem}
\usepackage{graphicx}
\usepackage{subfigure}
\usepackage{psfrag}
\newtheorem{thm}{Theorem}[section]

\newtheorem{cor}[thm]{Corollary}
\newtheorem{lem}[thm]{Lemma}
\theorembodyfont{\rmfamily}
\def\pf{\bigskip\noindent {\bf Proof.}~~}

\def\dfn#1{{\sl #1}}

\def\less{\backslash}

\def\pf{\bigskip\noindent {\emph{Proof.}}~~}

\newcounter{counter}
\def\pr#1{(\ref{#1})}

\begin{document}
\title{Hadwiger's conjecture  for graphs with forbidden holes}
\author{
Zi-Xia Song\thanks{Corresponding Author. Email: Zixia.Song@ucf.edu}\,\,  and Brian Thomas$^\dagger$ \thanks{Supported by the  UCF Research and Mentoring Program (RAMP) for undergraduate students.}\thanks{Current address: Department of Mathematics, University of Virginia, Charlottesville, VA 22904}\,\, 
\\
Department  of Mathematics\\
 University of Central Florida\\
Orlando, FL 32816, USA\\
}
\maketitle
\begin{abstract}

 Given a graph $G$, the \dfn{Hadwiger number} of $G$, denoted by  $h(G)$,  is the largest integer $k$ such that $G$ contains the complete graph $K_k$ as  a minor.  A \dfn{hole} in $G$ is an induced cycle of length at least four. 
  Hadwiger's Conjecture from 1943 states that for every graph $G$, $h(G)\ge \chi(G)$, where $\chi(G)$ denotes the chromatic number of $G$.  In this paper we establish more evidence for Hadwiger's conjecture by showing that if  a  graph $G$ with independence number $\alpha(G)\ge3$ has  no hole of length between $4$ and $2\alpha(G)-1$, then $h(G)\ge\chi(G)$. We also prove that if  a graph $G$ with independence number $\alpha(G)\ge2$  has  no hole of length between $4$ and $2\alpha(G)$, then $G$ contains an odd clique minor of size $\chi(G)$, that is,  such a graph $G$ satisfies the odd  Hadwiger's conjecture. 

\end{abstract}

\section{Introduction}

All graphs considered in this paper are
finite,  and have no loops or parallel edges.  We begin with some definitions.  
Let $G$ be a graph. The \dfn{complement} $\overline{G}$ of $G$ is the graph with vertex set  $V(G)$, such that two vertices are adjacent in $G$ if and only if they are non-adjacent in $\overline{G}$.   A \dfn{clique} in $G$ is a set of vertices all pairwise adjacent. A \dfn{stable set} in $G$ is a set of vertices all pairwise non-adjacent.  We use $\chi(G)$, $\omega(G)$, and  $\alpha(G)$ to denote the chromatic number, the clique number, and the independence number of $G$, respectively. 
Given a graph $H$, we say that $G$ is \dfn{$H$-free} if $G$ has no induced subgraph isomorphic to $H$. For a family  $\mathcal{F}$ of graphs, we say that $G$ is $\mathcal{F}$-free if $G$ is $F$-free for every  $F\in \mathcal{F}$.   A graph  $H$ is a \dfn{minor} of a graph $G$ if  $H$ can be
 obtained from a subgraph of $G$ by contracting edges.  
In those circumstances we also say that  $G$ has an $H$ \dfn{minor}.  \medskip

In 1943, Hadwiger~\cite{had} conjectured that for every integer $t\ge0$,  every graph either can be $t$-colored or has a $K_{t+1}$ minor.  
Hadwiger's conjecture is perhaps the most famous conjecture in graph theory, as pointed out by Paul Seymour in his recent survey~\cite{Seymour}. It  suggests a far reaching generalization of the Four Color Theorem~\cite{4ct1,4ct2,4ct3}. Hadwiger's conjecture is trivially true for $t\le2$, and reasonably easy for
$t=3$, as shown by Dirac~\cite{Dirac1952}. However, for $t\ge4$, Hadwiger's conjecture
implies the Four Color Theorem. 
Wagner~\cite{wagner} proved that the case $t=4$ of Hadwiger's conjecture is, in fact,
 equivalent to the Four Color Theorem, and the same was shown for $t=5$
by Robertson, Seymour, and Thomas~\cite{RST} in 1993.   Hadwiger's conjecture remains open for $t\ge6$.   As pointed out by Paul Seymour~\cite{Seymour} in his recent survey on Hadwiger's conjecture, proving that graphs with no $K_7$ minor are $6$-colourable is thus the first case of Hadwiger's  conjecture that is still open.  It  is not even known yet whether  every graph with no $K_7$ minor is $7$-colorable. 
 Kawarabayashi and Toft~\cite{kt} proved that every graph with no $K_7$ or $K_{4,\, 4}$ minor is $6$-colorable. 
 Jakobsen~\cite{Jakobsen1972, Jakobsen1983} proved that every graph with no $K_7^{=}$ minor is  $6$-colorable and  every graph with no $K_7^{-}$ minor is  $7$-colorable,  where for any integer $p>0$,
 $K_p^{-}$ (resp. $K_p^{=}$) denotes the graph obtained from $K_p$ by removing one edge (resp. two edges).   Recently  Albar and Gon\c calves~\cite{AG2015} proved that every graph with no $K_7$ minor is $8$-colorable and every graph with no $K_8$ minor is $10$-colorable. Their proof is computer-assisted.  In~\cite{roleksong2},  Rolek and the first author obtained  a short and computer-free proof of Albar and Gon\c calves' results and extended it to the next step by showing that every graph with no $K_9$ minor is $12$-colorable. Rolek and Song~\cite{roleksong2} also proved  that  every graph with no $K_8^=$ minor is $8$-colorable and every graph with no $K_8^-$ minor is $9$-colorable.   For more information on Hadwiger's conjecture, the readers are referred to an earlier  survey   by Toft~\cite{Toft}   and a very recent  informative survey  due to Seymour~\cite{Seymour}.    \medskip

 Hadwiger's conjecture has also been verified to be true for some special classes of graphs. In 2004, Reed and Seymour~\cite{rs2004} proved that Hadwiger's conjecture holds  for line graphs (possibly with parallel edges).  A graph $G$  is a \dfn{quasi-line graph} if for every vertex $v$, the  set of neighbors of $v$ can be covered by two  cliques, namely the vertex set of the neighborhood of $v$ can be partitioned into two cliques.  A graph $G$ is \dfn{claw-free} if $G$ does not contain $K_{1,3}$ as an induced subgraph. It is easy to verify that the class of line graphs is a proper subset of the class of quasi-line graphs, and the class of quasi-line graphs is a proper subset of the class of claw-free graphs. Recently quasi-line graphs attracted more attention, see~\cite{CP, CF2007, CF2008, ELTquasi}. In particular,  Chudnovsky and Seymour~\cite{CP} gave a constructive characterization of quasi-line graphs; Chudnovsky and Fradkin~\cite{CF2008, CF2010} proved in 2008 that Hadwiger's conjecture holds for quasi-line graphs and proved in 2010 that if $G$ is claw-free, then $G$ contains a clique minor of size at least $\lceil\frac{2}3\chi(G)\rceil$.  \medskip

One particular interesting case of Hadwiger's conjecture is  when graphs have independence number two. It has attracted more attention recently (see Section~4 in Seymour's survery~\cite{Seymour} for more information). As  stated in his survey,  Seymour believes that if Hadwiger's conjecture is true for graphs $G$ with $\alpha(G)=2$, then it is probably true  in general. Plummer, Stiebitz, and Toft~\cite{pst} proved that Hadwiger's conjecture holds for every $H$-free graph $G$ with $\alpha(G)=2$, where $H$ is any graph on four vertices and $\alpha(H)=2$. Later, Kriesell~\cite{kriesell} extended their result and proved that Hadwiger's conjecture holds for every $H$-free graph $G$ with $\alpha(G)=2$, where $H$ is any graph on five vertices and $\alpha(H)=2$.   \medskip

One strengthening of Hadwiger's conjecture is to consider the odd-minor variant of Hadwiger's conjecture. 
We say that a graph $G$ has an \dfn{odd clique minor} of size at least $k$  
if there are $k$ vertex-disjoint trees in $G$ such that every two of them are 
joined by an edge, and in addition, all the vertices of the trees are two-colored 
in such a way that the edges within the trees are bichromatic, but the edges 
between trees are monochromatic  (and hence the vertices of all trivial trees
 must receive the same color, where a tree is \dfn{trivial} if it has  one vertex only). 
We say that $G$ has an 
\dfn{odd $K_k$ minor} 
 if $G$ has an odd clique 
 minor of size at least $k$.  It is easy to see that any graph that  has an 
odd $K_k$ minor
certainly contains  $K_k$ as a minor. 
\medskip

Gerards and Seymour (see \cite{jt}, page 115.) proposed 
a well-known strengthening of  Hadwiger's  conjecture: for every integer $t\ge0$,  every graph either can be $t$-colored or has an odd $K_{t+1}$ minor.  This conjecture is referred to as ``the odd Hadwiger's conjecture". 
The odd Hadwiger's conjecture   is substantially stronger than Hadwiger's Conjecture. It    is trivially true for  $t\le 2$. 
The case $t=3$ was proved by  Catlin~\cite{catlin} in 1978.  
Guenin~\cite{guenin} announced  at a meeting in Oberwolfach in 2005  a solution of the case $t=4$. 
It  remains open for  $t\ge5$. Kawarabayashi and the first author~\cite{ks2} proved that  every graph $G$ on $n$ vertices with $\alpha(G)\ge2$ has an odd clique minor of size at least $\lceil n/(2\alpha(G)-1)\rceil$.   For more information on the odd Hadwiger's conjecture, the readers are referred to the very recent  survey  of  Seymour~\cite{Seymour}. \medskip

In this paper, we establish more evidence for Hadwiger's conjecture and the odd Hadwiger's conjecture.  
We first prove that if  a graph $G$ with independence number $\alpha(G)$  
has  no induced cycles of length between $4$ and $2\alpha(G)$, then $G$ satisfies the odd Hadwiger's conjecture. We then prove that if  a graph $G$ with independence number $\alpha(G)\ge3$  has  no induced cycles of length between $4$ and $2\alpha(G)-1$, then $G$ satisfies Hadwiger's conjecture.  We prove these two results in Section 2.\medskip

We need to introduce  more notation. Given a graph $G$, the \dfn{Hadwiger number} (resp. \dfn{odd Hadwiger number}) of $G$, denoted by  $h(G)$ (resp. $oh(G)$),  is the largest integer $k$ such that $G$ contains the complete graph $K_k$ as  a minor (resp. an odd minor). 
We use $|G|$, $\delta(G)$, and $\Delta(G)$ to denote the number
of vertices,  the  minimum degree, and the maximum degree of $G$, respectively. Given vertex sets $A, B \subseteq V(G)$, we say that $A$ is \dfn{complete to} (resp. \dfn{anti-complete to}) $B$ if for every $a \in A$ and every $b \in B$, $ab \in E(G)$ (resp. $ab \notin E(G)$).
The subgraph of $G$ induced by $A$, denoted $G[A]$, is the graph with vertex set $A$ and edge set $\{xy \in E(G) : x, y \in A\}$.  We denote by $G \less A$ the subgraph of $G$ induced on $V(G) \less A$.  If $A = \{a\}$, we simply write  $G \less a$.
An \dfn{$(A,B)$-path} of $G$ is a path with one end in $A$, the other end in $B$, and no internal vertices in $A\cup B$.  
A graph $G$ is \dfn{$t$-critical} if $\chi(G)=t$ and $\chi(G\less v)<t$ for any $v\in V(G)$. One can easily see that for any $t$-critical graph $G$, $\delta(G)\ge t-1$.     A \dfn{hole} in a graph $G$ is an induced cycle of length at least four.  An \dfn{antihole} in $G$ is an induced subgraph isomorphic to the  complement of  a hole.
A graph $G$ is \dfn{perfect} if $\chi(H)=\omega(H)$ for every induced subgraph $H$ of $G$.  We use $C_k$ to denote a cycle on $k$ vertices. We shall need the following two results. Theorem~\ref{spgt} is the  well-known Strong Perfect Graph Theorem~\cite{spgt}  and Theorem~\ref{quasi} is a result of Chudnovsky and Fradkin~\cite{CF2008}. 

\begin{thm}\label{spgt}
A graph $G$ is perfect if and only if it has no odd hole and no odd antihole.
\end{thm}
\begin{thm}\label{quasi}
If $G$ is a quasi-line graph, then $h(G)\ge\chi(G)$. 
\end{thm}

We shall need the following corollary. 
\begin{cor}\label{spgt1}
If $G$  is  $\{C_4, C_5, C_7, \dots,  C_{2\alpha(G)+1}\}$-free, then $G$ is perfect.
\end{cor}
\pf One can easily see that  $G$ is $C_k$-free for all $k \ge 2\alpha(G) +2$. Since $G$ is $C_4$-free, we see that $G$ is $\overline{C}_k$-free for all $k \ge 7$. Note that $C_5=\overline{C}_5$. Hence $G$ is $C_k$-free and $\overline{C}_k$-free for all $k \ge 5$  odd and so $G$ is perfect by Theorem~\ref{spgt}, as desired. 
\hfill\vrule height3pt width6pt depth2pt\\

\section{Main Results}

A graph $G$ is  an \dfn{inflation of a graph $H$}  if $G$ can be  obtained from $H$ by replacing each vertex  of $H$ by a clique of order at least one and two such cliques are complete to each other if their corresponding vertices in $H$ are adjacent. Under these  circumstances, we define $s(G)$ to be the size of the smallest clique used to replace a vertex of $H$. We first prove a lemma.

\begin{lem}\label{inflation}
Let $G$ be an inflation of an odd cycle $C$. Then   $\chi(G)\le \omega(G)+s(G)$ and $oh(G)\ge\chi(G)$.
\end{lem}
\pf Let $G$ and $C$ be  given as in the statement and let  $x_0, x_1, x_2, \dots, x_{2t}$ be the vertices of $C$ in order, where $t\ge1$ is an integer. The statement is trivially true when $t=1$. So we may assume that $t\ge2$.   Let $A_0, A_1, A_2, \dots, A_{2t}$ be vertex-disjoint cliques such that $A_i$ is used to replace $x_i$ for all $i\in\{0,1, \dots, 2t\}$. We may assume that $s(G)=|A_0|$. Then 
$|A_0|=\min\{|A_i|: \, 0\le i\le 2t\}$. Since $t\ge2$, we have  $\omega(G)=\omega(G\less A_0)$.  One can easily see that $G\less A_0$ is an inflation of a path on $2t$ vertices and so $\chi(G\less A_0)=\omega(G\less A_0)$. Therefore $\chi(G)\le\chi(G\less A_0)+|A_0|\le\omega(G\less A_0)+|A_0|=\omega(G)+s(G)$. \medskip

We next show that $G$ contains an odd clique minor of size $\omega(G)+|A_0|$. 
 Without loss of generality, we may assume that $\omega(G\less A_0)=|A_1|+|A_2|$. 
  By the choice of $|A_0|$, there are $|A_0|$ pairwise vertex-disjoint 
 $(A_0, A_3)$-paths (each  on $2t -1$ vertices) in $G\less (A_1\cup A_{2})$. Let $T_1, \dots ,T_{|A_0|}$ be  $|A_0|$   such paths and let $T_{|A_0|+1}, \dots,
   T_{|A_0|+\omega(G)}$ be $\omega(G)$ pairwise disjoint trivial trees (i.e., trees with one vertex only) in $A_1\cup A_2$. Now coloring  all the vertices in $A_{4} \cup A_{6} \cup \cdots \cup A_{2t}$ by color 1 and all the other vertices in $G$ by color 2, we see that  $T_1, \dots , T_{|A_0|}, T_{|A_0|+1}, \dots , T_{|A_0|+\omega(G)}$ yield an odd clique minor of size $|A_0|+\omega(G)\ge \chi(G)$, as desired. 
\hfill\vrule height3pt width6pt depth2pt\\
We next use Lemma~\ref{inflation} to prove that every graph $G$ with $\alpha(G)\ge2$ and no hole of length between $4$ and $2\alpha(G)$ satisfies the odd Hadwiger's conjecture, the following.\medskip

\begin{thm}\label{t1}
Let $G$ be a graph   with $\alpha(G) \ge 2$. 
If  $G$ is $\{C_4, \,C_5, \, C_6, \, \dots , \, C_{2\alpha(G)}\}$-free, 
 then $oh(G)\ge\chi(G)$. \end{thm}

\pf Let $\alpha = \alpha(G)$. We first assume that $G$ contains no odd hole of length $2\alpha +1$.
By Corollary~\ref{spgt}, $G$ is perfect and so  $G$ contains a clique (and thus an odd clique minor) of size $\chi(G)$.  So we may assume that $G$ contains an odd hole of length $2\alpha +1$,  say $C$, with vertices $ v_0, v_1, \dots, v_{2\alpha}$ in order.  We next prove that  for every  $w\in V(G\less C)$,  $w$  is either  complete to $C$ or adjacent to exactly three consecutive vertices on $C$.
Since  $\alpha(G)=\alpha$, we see that $w$ is adjacent to at least one vertex on $C$. If $w$ is complete to $C$, then we are done. So we may assume that  $wv_0\notin E(G)$ but $wv_1\in E(G)$. Then $w$ is not adjacent to $v_{2\alpha}, v_{2\alpha -1}, \dots , v_4$ as $G$ is $\{C_4, \,C_5, \, \dots , \, C_{2\alpha}\}$-free. Hence $w$ must be adjacent to $v_2, v_3$ because $\alpha(G)=\alpha$. This proves that  for any $w\in V(G\less C)$ not complete to $C$,  $w$ is adjacent to precisely three consecutive vertices on $C$.  \medskip

Let $J$ (possibly empty) denote the set of vertices in $G$ that are complete to $C$, and $A_i$ (possibly empty) denote the set of vertices in $G$ adjacent to $v_i, v_{i+1}, v_{i+2}$, where $i=0,1,\dots, 2\alpha$ and all arithmetic on indices here and henceforth is done modulo $2\alpha +1$.  Since $G$ is $C_4$-free, we see that $G[J]$ is a clique and 
$G[A_i]$ is a clique for all $i\in\{0,1,\dots, 2\alpha\}$. Note that $\{J, V(C), A_0, A_1, \dots , A_{2\alpha}\}$ partitions $V(G)$.  Since $\alpha(G) = \alpha$  and $G$ is $\{C_4, \,C_5, \, \dots , \, C_{2\alpha}\}$-free, one can easily check that $A_i$ is complete to $A_{i+1}\cup A_{i-1}$, and anti-complete to all $A_j$, where $j\ne i+1,i, i-1$. We now show  that $J$ is complete to all $A_i$. Suppose there exist $a\in J$ and $b\in A_i$ for some $i\in\{0, 1,\dots ,2\alpha\}$, say $i=0$,  such that $ab\notin E(G)$. Then $G[\{a, v_0, b, v_2\}]$ is an induced $C_4$, a contradiction.  Thus $J$ is complete to all $A_i$ and so $J$ is complete to $V(G)\less J$.  Let $A_i^*=A_i\cup \{v_{i+1}\}$ for $i=0, 1, \dots ,2\alpha$. Then $A_i^*\not= \emptyset$ for all $i\in\{0,1,\dots, 2\alpha\}$  and  $G\less J$ is 
 an inflation of an odd cycle of length $2\alpha +1$.  By Lemma~\ref{inflation}, $oh(G\less J)\ge \chi(G\less J)$ and so $oh(G)\ge oh(G\less J)+|J|\ge\chi(G\less J)+|J|=\chi(G)$, as desired. \medskip
 
 This completes the proof of Theorem~\ref{t1}.  \hfill\vrule height3pt width6pt depth2pt\\

Finally we prove that every graph $G$ with $\alpha(G)\ge3$ and no hole of length between $4$ and $2\alpha(G)-1$ satisfies Hadwiger's conjecture.\medskip

\begin{thm}\label{main}
Let $G$ be a graph  with  $\alpha(G)\ge 3$. If $G$ is $\{C_4, \, C_5, C_6, \dots, C_{2\alpha(G)-1}\}$-free, then $h(G)\ge  \chi(G)$. 
\end{thm}

\setcounter{counter}{0}

\pf   Suppose for a contradiction that $h(G)<  \chi(G)$.   Let  $G$ be a counterexample with $|V(G)|$ as small as possible.  Let $n :=|V(G)|$,  $t :=\chi(G)$,   and  $\alpha :=\alpha(G)$.  By Theorem~\ref{spgt} and the fact that $h(G)\ge\omega(G)$, $G$ is not perfect. Since  $G$ is $\{C_4, \, C_5, C_7, \dots, C_{2\alpha-1}\}$-free, by Corollary~\ref{spgt1},   we see that   $G$ must contain an odd hole, say $C$, with  $2\alpha+1$ vertices. \medskip

\noindent \refstepcounter{counter}\label{e:mindeg} (\arabic{counter}) $G$ is $t$-critical and so $\delta(G)\ge t-1$.

\pf Suppose that there exists  $x\in V(G)$ such that  $\chi(G\less x)=t$. If $\alpha(G\less x)=\alpha$, then $G\less x$ is $\{C_4, \, C_5, C_6, \dots, C_{2\alpha(G\less x)-1}\}$-free. By the minimality of $G$,  we have $h(G\less x)\ge \chi(G\less x)=t$ and so $h(G)\ge h(G\less x)\ge t$, a contradiction. Thus $\alpha(G\less x)=\alpha-1$. Then $G\less x$ is $\{C_4, \, C_5, C_7, \dots, C_{2\alpha(G\less x)+1}\}$-free. By Corollary~\ref{spgt1}, $G\less x$ is perfect and so $h(G)\ge h(G\less x)\ge\omega(G\less x)=\chi(G\less x)=t$, a contradiction. Thus $G$ is $t$-critical and so  $\delta(G)\ge t-1$. \hfill\vrule height3pt width6pt depth2pt\medskip

\noindent \refstepcounter{counter}\label{e:maxdeg} (\arabic{counter})  $\Delta(G)\le n-2$.

\pf Suppose there exists a vertex $x$ in $G$ with $d(x)=n-1$. Then $\chi(G\less x)=t-1$ and by the minimality of $G$, $h(G\less x)\ge \chi(G\less x)=t-1$ and so $h(G)\ge h(G\less x) +1\ge (t-1)+1=t$,  a contradiction. \hfill\vrule height3pt width6pt depth2pt\medskip

\noindent \refstepcounter{counter}\label{e:clique} (\arabic{counter}) $\omega(G)\le t-2$.

\pf Suppose that $\omega(G)\ge t-1$. Since  $h(G)< t$, we see that  $ \omega(G)=t-1$.  Let $H\subseteq G$ be isomorphic to $K_{t-1}$. Then $C$ and $H$ have at most two vertices in common, and if $|C\cap H|=2$, then the two vertices in $C\cap H$ must be adjacent on $C$. Let $P$ be a subpath of $C\less H$ on $2\alpha-1$ vertices. Then $P$ is an induced path in  $G\less H$.  Since $\alpha(G)=\alpha$, we see that every vertex in $H$ must have a neighbor on $P$.  By contracting the path $P$ into a single vertex, we see that  $h(G)\ge t$, a contradiction.
\hfill\vrule height3pt width6pt depth2pt\\

Let  $ v_0, v_1, \dots, v_{2\alpha}$ be the   vertices of $C$ in order.  We next show that \medskip

\noindent \refstepcounter{counter}\label{e:nbr} (\arabic{counter})  For every  $w\in V(G\less C)$,  either $w$ is complete to $C$, or $w$ is adjacent to exactly three consecutive vertices on $C$, or $w$ is adjacent to exactly four consecutive vertices on $C$. 

\pf Since  $\alpha(G)=\alpha$, we see that $w$ is adjacent to at least one vertex on $C$. Suppose that $w$ is not complete to $C$. We may assume that  $wv_0\notin E(G)$ but $wv_1\in E(G)$. Then $w$ is not adjacent to $v_{2\alpha}, v_{2\alpha -1}, \dots , v_5$ because  $G$ is $\{C_4, \,C_5, \, \dots , \, C_{2\alpha-1}\}$-free.    If $wv_4\in E(G)$, then $w$ must be adjacent to $v_2, v_3$ because  $G$ is $\{C_4, \, C_5\}$-free. If $wv_4\notin E(G)$, then again $w$ must be adjacent to $v_2, v_3$  because $\alpha(G)=\alpha$.  Thus $w$ is adjacent to either  $v_1, v_2, v_3$  or $v_1, v_2,v_3, v_4$  on $C$,  as desired. \hfill\vrule height3pt width6pt depth2pt\medskip

Let $J$  denote the set of vertices in $G$ that are complete to $C$.  For each $i\in I : =\{0,1,\dots, 2\alpha\}$, let  $A_i\subseteq V(G\less C)$ (possibly empty) denote the set of vertices in $G$ adjacent to precisely $v_i, v_{i+1}, v_{i+2}$ on $C$, and let $B_i\subseteq V(G\less C)$ (possibly empty) denote the set of vertices in $G$ adjacent to precisely $v_i, v_{i+1}, v_{i+2}, v_{i+3}$ on $C$, where all arithmetic on indices here and henceforth is done modulo $2\alpha+1$. By \pr{e:nbr}, $\{J, V(C), A_0, A_1, \dots, A_{2\alpha}, B_0, B_1,  \dots, B_{2\alpha}\}$ partitions $V(G)$. \medskip

\noindent \refstepcounter{counter}\label{e:J} (\arabic{counter})
   $J=\emptyset$.

 \pf Suppose that $J\not=\emptyset$.  Let   $a\in J$.  By \pr{e:maxdeg}, there exists $b\in V(G)\less (V(C)\cup J)$  such that $ab\notin E(G)$. By \pr{e:nbr}, we may assume that $b$ is adjacent to $v_0, v_1, v_2$. Then $G[\{a, v_0, b, v_2\}]$ is an induced $C_4$ in $G$, a contradiction. \hfill\vrule height3pt width6pt depth2pt\medskip

  The fact that $\alpha(G) = \alpha$ implies that \\

\noindent \refstepcounter{counter}\label{e:ABclique} (\arabic{counter}) For each $i\in I$, both $G[A_i]$ and $G[B_i]$ are cliques; and $A_i$ is complete to $A_{i-1}\cup A_{i+1}$.\\

Since $G$ is $\{C_4, \,C_5, \, \dots , \, C_{2\alpha-1}\}$-free, one can easily check that \\

\noindent \refstepcounter{counter}\label{e:Aac} (\arabic{counter})
For each $i\in I$, $A_i$ is  anti-complete to each $A_j$, where $j\in I\less\{i-2, i-1, i, i+1, i+2\}$; and \\

\noindent \refstepcounter{counter}\label{e:Bac} (\arabic{counter})
For each $i\in I$, $B_i$ is complete to $B_{i-1}\cup A_{i}\cup A_{i+1}\cup B_{i+1}$ and anti-complete to each $B_j$, where $j\in I\less\{ i-1, i, i+1\}$.\\

We shall also  need the following:\\

\noindent \refstepcounter{counter}\label{e:oneneighbor} (\arabic{counter})
For each $i\in I$, if $B_i\not=\emptyset$, then  $B_{j}=\emptyset$ for any $j\in I\less\{i-2, i-1, i, i+1, i+2\}$.

\pf    Suppose $B_j$ is not empty for some $j\in I\less\{ i+2, i+1, i, i-1, i-2\}$.  We may assume that $j>i$. Let $a\in B_i$ and $b\in B_{j}$. Then $G[\{v_i, a, v_{i+3},\dots, v_j,  b, v_{j+3}, \dots, v_{i-1}\}]$ is an odd hole with $2\alpha-1$ vertices  if $ab\notin E(G)$, and $G[\{v_i, a,   b, v_{j+3}, \dots, v_{i-1}\}]$  is a hole  with  $2\alpha+1-j+i\ge4$ vertices  if $ab\in E(G)$. In both cases, we obtain a contradiction.  \hfill\vrule height3pt width6pt depth2pt\\

With an argument similar to that  of  \pr{e:oneneighbor}, we see that  \\

\noindent \refstepcounter{counter}\label{e:Nneighbor} (\arabic{counter})
For each $i\in I$, if $B_i\not=\emptyset$, then  $B_{i}$ is anti-complete to $A_j$ for any  $j\in I\less\{i-1, i, i+1, i+2\}$.\\

We next show that \\

\noindent \refstepcounter{counter}\label{e:partition} (\arabic{counter})
For any $i\in I$, if $A_i\ne\emptyset$, then each vertex in $A_i$ is either anti-complete to $A_{i+2}$ or anti-complete to $A_{i-2}$.

\pf Suppose there exists a vertex $x\in A_i$ such that  $x$ is adjacent to a vertex $y\in A_{i-2}$ and a vertex $z\in A_{i+2}$. Then $G[ \{x, y,z\}\cup (V(C)\less \{v_{i-1}, v_i, v_{i+1}, v_{i+2}, v_{i+3}\})]$ is  an odd hole with  $2\alpha-1$ vertices, a contradiction. \hfill\vrule height3pt width6pt depth2pt\\

\noindent \refstepcounter{counter}\label{e:A{j-1}A{j+2}} (\arabic{counter})
For any $i\in I$,    if $B_i\ne\emptyset$,  then every vertex in $B_i$ is either complete to $A_{i-1}$ or complete to  $A_{i+2}$.

\pf Suppose for a contradiction,  say $B_2\ne\emptyset$,  and  there exists a vertex $b\in B_2$  such that $b$ is not adjacent to a vertex $a_1\in A_1$ and a vertex $a_4\in A_4$. By \pr{e:Aac}, $A_1$ is anti-complete to $A_4$. Thus $G$ contains a stable set $\{b, a_1, a_4,v_0\}$ of size four if $\alpha=3$ or stable  set $\{b, a_1, a_4, v_0, v_7, v_9, \dots, v_{2\alpha-1}\}$ of size $\alpha+1$ if $\alpha\ge4$, a contradiction. \hfill\vrule height3pt width6pt depth2pt\\

\noindent \refstepcounter{counter}\label{e:3B} (\arabic{counter})
There exists an $i\in I$ such that $B_j= \emptyset$ for any $j\in I\less \{i, i+1, i+2\}$.

\pf This is obvious  if $B_k=\emptyset$ for any $k\in I$. So we may assume that $B_k\ne\emptyset$ for some $k\in I$, say $B_2\neq \emptyset$.  Then  by \pr{e:oneneighbor},  $B_j=\emptyset$  for all  $j=5, 6, \dots, 2\alpha$.    By \pr{e:oneneighbor} again, either $B_0\ne \emptyset$ or $B_4\ne \emptyset$ but not both.  By symmetry,  we may assume  that $B_4=\emptyset$. Similarly, either $B_0\ne \emptyset$ or $B_3\ne \emptyset$ but not both. Thus either $B_j=\emptyset$ for all $j\in I\less \{0,1,2\}$ or  $B_j=\emptyset$ for all $j\in I\less \{1,2, 3\}$. \hfill\vrule height3pt width6pt depth2pt\\

By \pr{e:3B}, we may assume that  $B_j=\emptyset$ for all $j\in I\less \{1,2,3\}$. 
For any  $A_i\ne\emptyset$, where $i\in I$, 
 let $A_i^1=\{a\in A_i: a  \text{  has a neighbor in } A_{i-2}\}$, $A_i^3=\{a\in A_i: a\, \text{ has a neighbor in } A_{i+2}\}$,  and $A_i^2=A_i\less (A_i^1\cup A_i^3)$. Then $A_i^2$ is anti-complete to $A_{i-2}\cup A_{i+2}$. By \pr{e:partition}, 
  $A_i^1$ is anti-complete to $A_{i+2}$ and $A_i^3$ is anti-complete to $A_{i-2}$.   Clearly, $\{A_i^1, A_i^2, A_i^3\}$ partitions $A_i$. 
  Next,  for any  $B_j\ne\emptyset$, where $j\in \{1,2,3\}$, by \pr{e:A{j-1}A{j+2}}, let $B_j^1=\{b\in B_j: b\, \text{ is complete to }  A_{j-1}\}$ and  $B_j^2=\{b\in B_j: b\, \text{ is complete to }  A_{j+2}\}$. Clearly,  $B_j^1$ and $B_j^2$ are not necessarily disjoint. Note that  $B_j^1$ and $B_j^2$ are not symmetrical because  $B_j^1$  is complete to $A_{j-1}$ and $B_j^2$  is complete to $A_{j+2}$.  \\
   
 \noindent \refstepcounter{counter}\label{e:BjA{j-1}} (\arabic{counter})
 For any $j\in\{1,2,3\}$, $B_j$ is anti-complete to $A_{j-1}^1\cup A_{j+2}^3$.  
 
 \pf Suppose  there exist a vertex $b\in B_j$ and a vertex $a\in A_{j-1}^1\cup A_{j+2}^3$ such that $ba\in E(G)$. 
 By the definitions of $A_{j-1}^1$ and $A_{j+2}^3$, we see that $a$ has a neighbor, say $c$, in $A_{j-3}$ if $a\in A_{j-1}^1$, or  in $A_{j+4}$ if $a\in A_{j+2}^3$. Now $G[ \{b, a,c\}\cup (V(C)\less \{v_{j-2}, v_{j-1}, v_j, v_{j+1}, v_{j+2}\})]$ is  an odd hole  of length $2\alpha-1$ if $a\in A_{j-1}^1$, or  $G[ \{b, a,c\}\cup (V(C)\less \{v_{j+1}, v_{j+2}, v_{j+3}, v_{j+4}, v_{j+5}\})]$  is  an odd hole of length $2\alpha-1$ if $a\in A_{j+2}^3$. In either case, we obtain a contradiction. \hfill\vrule height3pt width6pt depth2pt\\

 \noindent \refstepcounter{counter}\label{e:quasi-line} (\arabic{counter})
 $G$ is a quasi-line graph.
 
 \pf   It suffices to show that for any $x\in V(G)$, $N(x)$ is covered by two  cliques.   By \pr{e:J}, $J=\emptyset$. Since $B_j^1$ and $B_j^2$ are not symmetrical  for all $j\in\{1,2,3\}$, we consider the following four cases. \\
 
\noindent{\bf Case 1:}   $x\in A_i$ for some $i\in I\less\{0, 1,2,3,4,5\}$.\\

In this case, $x\in A_i^k$ for some $k\in\{1,2,3\}$. We first assume that $k=1$. Then $x\in A_i^1$. By \pr{e:partition} and the definition of $A_i^1$,  $x$ is anti-complete to $A_{i+2}$ and so  $N[x]\subseteq A_{i-2}\cup A_{i-1}\cup A_i\cup A_{i+1}\cup\{v_{i}, v_{i+1}, v_{i+2}\}$. We see that $N(x)$ is covered by two  cliques $G[(A_{i-2}\cap N(x))\cup A_{i-1}\cup \{v_i\}]$ and $G[(A_i\less x)\cup A_{i+1}\cup\{v_{i+1}, v_{i+2}\}]$. By symmetry, the same holds if $k=3$. So we may assume that $k=2$. By the definition of $A_i^2$, $x$ is anti-complete to $A_{i-2}\cup A_{i+2}$. Thus  $N[x]=A_{i-1}\cup A_i\cup A_{i+1}\cup\{v_i, v_{i+1}, v_{i+2}\}$ and so $N(x)$ is covered by two  cliques  $G[A_{i-1}\cup (A_i\less x)\cup\{v_i\}]$  and $G[A_{i+1}\cup \{v_{i+1}, v_{i+2}\}]$.   \\

\noindent{\bf Case 2:}   $x\in A_i$ for some $i\in \{0,1,2,3,4,5\}$.\\

In this case,  we first assume that $i=0$. Then $x\in A_0^k$ for some $k\in\{1,2,3\}$. Assume that  $x\in A_0^1$. Then  $x$ is anti-complete to $B_1$ by \pr{e:BjA{j-1}} and anti-complete to $A_2$ by \pr{e:partition}. Thus  $N[x]\subseteq A_{2\alpha-1}\cup A_{2\alpha}\cup A_0\cup A_1\cup \{v_0, v_1, v_2\}$.  One can see that $N(x)$ is covered by two  cliques   $G[(A_{2\alpha-1}\cap N(x))\cup A_{2\alpha}\cup \{v_0\}]$ and $G[(A_0\less x)\cup A_1\cup\{v_1, v_2\}]$.  It can be easily checked that $N(x)$ is covered by two  cliques $G[A_{2\alpha}\cup (A_0\less x)\cup\{v_0\}]$ and $G[A_1\cup (B_1\cap N(x))\cup\{v_1, v_2\}]$ if $k=2$; and  by two  cliques  $G[A_{2\alpha}\cup (A_0\less x)\cup\{v_0\}]$ and $G[A_1\cup (B_1\cap N(x))\cup(A_2\cap N(x))\cup\{v_1, v_2\}]$ if $k=3$.  \\

Next assume that $i=1$. Then $x\in A_1^k$ for some $k\in\{1,2,3\}$ and $N[x]\subseteq A_{2\alpha}\cup A_0\cup A_1\cup B_1\cup A_2\cup B_2\cup A_3\cup\{v_1, v_2, v_3\}$.  One can see that $N(x)$ is covered by two  cliques   $G[(A_{2\alpha}\cap N(x))\cup A_0\cup \{v_1\}]$ and $G[(A_1\less x)\cup B_1\cup A_2\cup\{v_2, v_3\}]$ if $k=1$; by two  cliques  $G[A_0\cup (A_1\less x)\cup \{v_1\}]$ and $G[B_1\cup A_2\cup (B_2\cap N(x))\cup\{v_2, v_3\}]$ if $k=2$; and two  cliques  $G[A_0\cup (A_1\less x)\cup B_1^1\cup\{v_1, v_2\}]$ and $G[B_1^2\cup A_2\cup (B_2\cap N(x))\cup (A_3\cap N(x))\cup\{ v_3\}]$ if $k=3$. \medskip

Assume that $i=2$. Then $x\in A_2^k$ for some $k\in\{1,2,3\}$ and $N[x]\subseteq A_{0}\cup A_{1}\cup B_{1}\cup A_2\cup B_2\cup A_{3}\cup B_{3}\cup A_4\cup\{ v_2, v_3, v_4\}$.  One can check  that $N(x)$ is covered by two  cliques   $G[(A_0\cap N(x))\cup A_1\cup B_1^1\cup \{v_2\}]$ and $G[(A_2\less x)\cup B_1^2\cup B_2\cup A_3\cup\{v_3, v_4\}]$ if $k=1$; by two  cliques  $G[A_1\cup B_1\cup (A_2\less x)\cup \{v_2\}]$ and $G[B_2\cup A_3\cup (B_3\cap N(x))\cup\{v_3, v_4\}]$ if $k=2$;  and by two cliques  $G[A_1\cup B_1\cup (A_2\less x)\cup B_2^1\cup\{v_2, v_3\}]$ and $G[B_2^2\cup A_3\cup (B_3\cap N(x))\cup (A_4\cap N(x))\cup\{ v_4\}]$ if $k=3$.  \\

Assume that $i=3$. Then $x\in A_3^k$ for some $k\in\{1,2,3\}$ and $N[x]\subseteq A_1\cup B_{1}\cup A_2\cup B_2\cup A_{3}\cup B_{3}\cup A_4\cup A_5\cup\{ v_3, v_4, v_5\}$.  One can check  that $N(x)$ is covered by two  cliques   $G[(A_1\cap N(x)) \cup (B_1\cap N(x))\cup A_2\cup B_2^1\cup \{v_3\}]$ and $G[B_2^2\cup (A_3\less x)\cup   B_3\cup A_4\cup\{v_4, v_5\}]$ if $k=1$; and by two  cliques  $G[(B_1\cap N(x))\cup A_2\cup B_2\cup \{v_3\}]$ and $G[(A_3\less x)\cup B_3 \cup A_4\cup\{v_4, v_5\}]$ if $k=2$. So we may assume that $x\in A_3^3$. By \pr{e:BjA{j-1}}, we have $A_3^3$ is anti-complete to $B_1$. Thus $N(x)$ is covered by  two  cliques  $G[A_2\cup B_2\cup (A_3\less x)\cup B_3^1\cup\{v_3, v_4\}]$ and $G[B_3^2\cup A_4\cup (A_5 \cap N(x))\cup\{ v_5\}]$ if $k=3$.  \\

Assume that $i=4$. Then $x\in A_4^k$ for some $k\in\{1,2,3\}$ and $N[x]\subseteq A_2\cup B_{2}\cup A_3\cup B_3\cup  A_4\cup A_5\cup A_6\cup\{ v_4, v_5, v_6\}$.  One can see  that $N(x)$ is covered by two  cliques   $G[(A_2\cap N(x)) \cup (B_2\cap N(x))\cup A_3\cup B_3^1\cup \{v_4\}]$ and $G[B_3^2\cup (A_4\less x)\cup  A_5\cup\{v_5, v_6\}]$ if $k=1$, and  by two  cliques  $G[(B_2\cap N(x))\cup A_3\cup B_3\cup \{v_4\}]$ and $G[(A_4\less x) \cup A_5\cup\{v_5, v_6\}]$ if $k=2$. So we may assume that $x\in A_4^3$. By  \pr{e:BjA{j-1}}, we have $A_4^3$ is anti-complete to $B_2$.  Thus $N(x)$ is covered  by  two  cliques  $G[A_3\cup B_3\cup (A_4\less x)\cup\{v_4, v_5\}]$ and $G[ A_5\cup (A_6 \cap N(x))\cup\{ v_6\}]$ if $k=3$.  \\

Finally assume that $i=5$. Then $x\in A_5^k$ for some $k\in\{1,2,3\}$ and $N[x]\subseteq A_3\cup B_{3}\cup A_4\cup   A_5\cup A_6\cup A_7\cup\{ v_5, v_6, v_7\}$.  One can check  that $N(x)$ is covered by two  cliques   $G[(A_3\cap N(x)) \cup (B_3\cap N(x))\cup A_4\cup  \{v_5\}]$ and $G[(A_5\less x)\cup  A_6\cup\{v_6, v_7\}]$ if $k=1$,  and  by two  cliques  $G[(B_3\cap N(x))\cup A_4\cup \{v_5\}]$ and $G[(A_5\less x) \cup A_6\cup\{v_6, v_7\}]$ if $k=2$.  So we may assume that $x\in A_5^3$. By  \pr{e:BjA{j-1}}, we have $A_5^3$ is anti-complete to $B_3$.   Thus $N(x)$ is covered by  two  cliques  $G[ A_4\cup (A_5\less x)\cup \{v_5, v_6\}]$ and $G[ A_6\cup (A_7\cap N(x))\cup\{v_7\}]$  if $k=3$. \medskip

This completes the proof of Case 2. \\

\noindent{\bf Case 3:}   $x\in B_j$ for some $j\in \{1,2,3\}$.\\

In this case, first assume that $j=1$.  Then $x\in B_1^k$ for some $k\in\{1,2\}$, and  $N[x]\subseteq A_0\cup A_1\cup B_1\cup A_2\cup B_2\cup A_3\cup\{v_1, v_2, v_3, v_4\}$. We see that  $N(x)$ is covered by two  cliques   $G[A_0\cup A_1\cup (B_1^1\less x)\cup \{v_1, v_2\}]$ and $G[B_1^2\cup A_2\cup B_2\cup (A_3\cap N(x))\cup\{v_3, v_4\}]$ if $k=1$; and by two cliques  $G[(A_0\cap N(x))\cup A_1\cup B_1^1\cup \{v_1, v_2\}]$ and $G[(B_1^2\less x)\cup A_2\cup B_2\cup A_3\cup\{v_3, v_4\}]$ if $k=2$. \medskip

Next assume that $j=2$. Then $x\in B_2^k$ for some $k\in\{1,2\}$, and  $N[x]\subseteq A_1\cup B_1\cup A_2\cup B_2\cup A_3\cup B_3\cup A_4\cup \{ v_2, v_3, v_4, v_5\}$. One can  see that $N(x)$ is covered by two cliques   $G[A_1\cup B_1\cup A_2\cup (B_2^1\less x)\cup \{v_2, v_3\}]$ and $G[B_2^2\cup A_3\cup B_3\cup (A_4\cap N(x))\cup\{v_4, v_5\}]$ if $k=1$; and by two  cliques $G[(A_1\cap N(x))\cup B_1\cup A_2\cup B_2^1\cup \{v_2, v_3\}]$ and $G[(B_2^2\less x)\cup A_3\cup B_3\cup A_4\cup\{v_4, v_5\}]$ if $k=2$. \medskip

Finally assume $j=3$. Then $x\in B_3^k$ for some $k\in\{1,2\}$, and  $N[x]\subseteq A_2\cup B_2\cup A_3\cup B_3\cup A_4\cup A_5\cup \{v_3, v_4, v_5, v_6\}$. We see that $N(x)$ is covered by two cliques   $G[A_2\cup B_2\cup A_3\cup (B_3^1\less x)\cup \{v_3, v_4\}]$ and $G[B_3^2\cup A_4\cup (A_5\cap N(x))\cup\{v_5, v_6\}]$ if $k=1$; and by two  cliques  $G[(A_2\cap N(x))\cup B_2\cup A_3\cup B_3^1\cup \{v_3, v_4\}]$ and $G[(B_3^2\less x)\cup A_4\cup A_5\cup\{v_5, v_6\}]$ if $k=2$. \\

\noindent{\bf Case 4:}   $x\in V(C)$. \\

In this case, let $x=v_i$ for some $i\in I$.   First assume that $i\ne 1,2,3,4,5,6$.  Then  $N(v_{i})=A_{i-2}\cup A_{i-1}\cup A_{i}\cup\{v_{i-1}, v_{i+1}\}$  and so   $N(v_{i})$ is  covered by two   cliques $G[A_{i-2}\cup A_{i-1}\cup\{v_{i-1}\}]$ and $G[A_{i}\cup\{ v_{i+1}\}]$. Next assume that $i\in\{ 1,2,3,4,5,6\}$.   One can easily check that $N(v_1)$ is covered by two  cliques $G[A_{2\alpha}\cup A_0\cup\{v_0\}]$ and $G[A_1\cup B_1\cup\{v_2\}]$;  $N(v_2)$   by  two  cliques  $G[A_{0}\cup A_1\cup\{v_1\}]$ and $G[B_1\cup A_2\cup B_2\cup\{v_3\}]$;      $N(v_3)$  by two  cliques $G[A_{1}\cup B_1\cup A_2\cup\{v_2\}]$ and $G[B_2\cup A_3\cup B_3\cup\{v_4\}]$;  $N(v_4)$  by two  cliques $G[ B_1\cup A_2\cup B_2\cup\{v_3\}]$ and $G[A_3\cup B_3\cup A_4\cup\{v_5\}]$;  $N(v_5)$  by two  cliques $G[ B_2\cup A_3\cup B_3\cup\{v_4\}]$ and $G[A_4\cup  A_5\cup\{v_6\}]$;  and $N(v_6)$  by two  cliques $G[ B_3\cup A_4\cup\{v_5\}]$ and $G[A_5\cup  A_6\cup\{v_7\}]$, respectively. 
\medskip

This proves that $G$ is a quasi-line graph. \hfill\vrule height3pt width6pt depth2pt\\

By \pr{e:quasi-line}, $G$ is a quasi-line graph.  By Theorem~\ref{quasi}, $h(G)\ge \chi(G)=t$, a contradiction. This completes the proof of Theorem~\ref{main}.
  \hfill\vrule height3pt width6pt depth2pt\\

\noindent {\bf Remark.}  We made no use of the fact $\omega(G)\le t-2$ in the proof of Theorem~\ref{main}. We kept it in the proof in the hope that  one might be able to find a short proof of  Theorem~\ref{main} without using the fact that Hadwiger's conjecture is true for quasi-line graphs (namely Theorem~\ref{quasi}). 

\section*{Acknowledgements}
The authors would like to thank the anonymous referee for many helpful comments.

\end{document}